\documentclass[english,12pt,leqno]{article}
\usepackage{tikz}
\usetikzlibrary{decorations.markings}
\usetikzlibrary{calc}
\usepackage[latin1]{inputenc}
\usepackage{amsxtra}
\usepackage{amsmath}
\usepackage{amssymb}
\usepackage{amsfonts}
\usepackage[all]{xy}
\usepackage{mathrsfs}
\usepackage{amsthm}
\usepackage{mathabx}
\usepackage{float}

\usepackage{babel}

\newtheorem{thm}[equation]{Theorem}

\newtheorem{cor}[equation]{Corollary}
\newtheorem{lem}[equation]{Lemma}
\newtheorem{prop}[equation]{Proposition}

\newtheoremstyle{example}{\topsep}{\topsep}%
     {}
     {}
     {\bfseries}
     {.}
     {2pt}
     {\thmname{#1}\thmnumber{ #2}\thmnote{ #3}}

   \theoremstyle{example}
   
   \newtheorem{Defi}[equation]{Definition}
   
   \newtheorem{rem}[equation]{Remark}
   \newtheorem{rems}[equation]{Remarks}

   \newtheorem{ex}[equation]{Example}
   
    \newtheorem{refo}[equation]{Reformulation}

\newtheoremstyle{example}{\topsep}{\topsep}%
     {}
     {}
     {\bfseries}
     {.}
     {2pt}
     {\thmname{#1}\thmnumber{ #2}\thmnote{ #3}}

   \numberwithin{equation}{section}

\def\eps{{\varepsilon}}

\def\CC{\mathbb{C}}
\def\DD{\mathbb{D}}
\def\EE{\mathbb{E}}

\def\GG{\mathbb{G}}
\def\LL{\mathbb{L}}

\def\RR{\mathbb{R}}
\def\ZZ{\mathbb{Z}}

\def\HH{\mathbb{H}}

\def\Fen{\mathfrak{F}}

\def\hen{\mathfrak{h}}

\def\Sen{\mathfrak{S}}

 \def\Gen{\mathfrak{G}}

\def\Ac{\mathcal{A}}
\def\Bc{\mathcal{B}}
\def\Cc{\mathcal{C}}
\def\Kc{\mathcal{K}}
\def\Dc{\mathcal{D}}
\def\Ec{\mathcal{E}}
\def\Fc{\mathcal{F}}
\def\Gc{\mathcal{G}}
\def\Ic{\mathcal{I}}
\def\Jc{\mathcal{J}}
\def\Lc{\mathcal{L}}
\def\Mc{\mathcal{M}}
\def\Nc{\mathcal{N}}
\def\Hc{\mathcal{H}}
\def\Oc{\mathcal{O}}
\def\Pc{\mathcal{P}}
\def\Qc{\mathcal{Q}}
\def\Rc{\mathcal{R}}
\def\Sc{\mathcal{S}}

\def\Xc{\mathcal{X}}

\def\be{\begin{equation}}

\def\cod{\on{codim}}

\def\del{{\partial}}
\def\dim{{\rm{dim}}}

\def\DSh{D^b\Sh }

\def\ee{\end{equation}}

\def\Hom{\on{Hom}}
\def\hra{\hookrightarrow}

\def\Id{{\on{Id}}}

\def\k {\mathbf k}

\def\lra{\longrightarrow}

\def\ol{\overline}
\def\on{\operatorname}
\def\op{{\on{op}}}
\def\orr{\on{or}}
\def\OR{\on{or}}

\def\para{{\parallel}}
\def\Perv{\on{Perv}}

\def\Rep{\on{Rep}}

\def\sgn{\on{sgn}}
\def\Sh{\on{Sh}}
\def\Supp{{\on{Supp}}}

\def\ul{\underline}

\def\Vect{{\on{Vect}}}
\def\Vfd{\on{Vect}^{\on{fd}}_\k}

\def\wc{\widecheck}
  \def\wt{\widetilde}

\title{ Perverse sheaves over real hyperplane arrangements II
}

\author{{Mikhail Kapranov} \and {Vadim Schechtman}}

\begin{document}

\maketitle

\begin{abstract}
 Let $\Hc$ be an arrangement of hyperplanes in $\RR^n$ and $\Perv(\CC^n,\Hc)$ be
 the category of perverse sheaves on $\CC^n$ smooth with respect to the stratification
 given by complexified flats of $\Hc$. We give a description of $\Perv(\CC^n, \Hc)$ in
 terms of ``matrix diagrams'', i.e., diagrams formed by vector spaces $E_{A,B}$ labelled
 by pairs $(A,B)$ of real faces of $\Hc$ (of all dimensions) or, equivalently, by the cells
 $iA+B$ of a natural cell decomposition of $\CC^n$. A matrix diagram is formally similar to
 a datum describing a constructible (non-perverse) sheaf but with the direction of one half of
 the arrows reversed. 
\end{abstract}

\tableofcontents

 \addtocounter{section}{-1}
 

\section{Introduction}
\paragraph*{}
Let  $\Hc$ be a finite arrangement of linear hyperplanes in $\RR^n$. 
It gives rise to a stratification $\Sc^{(0)}=\Sc^{(0)}_\Hc$ of the complex
space $\CC^n$ into the generic parts of the complexified flats of $\Hc$,
see \S  \ref{sec:gen-arr} for a precise definition. 
  Let   $\k$ be a field.
We denote by $\Perv(\CC^n,\Hc)$  the category of perverse sheaves (middle perversity)
 of $\k$-vector spaces on $\CC^n$ which are constructible with respect to $\Sc^{(0)}$. 
 Such perverse sheaves and their categories are of great importance in several areas,
 including  representation theory of quantum groups \cite{BFS}.
 They also provide a large class of nontrivial examples of categories of perverse sheaves.

\vskip .2cm

In \cite{KS}, we gave a description of $\Perv(\CC^n,\Hc)$ 
in terms of certain quivers,
i.e., diagrams of vector spaces $E_A$ labelled by {\em faces}  $A$ of $\Hc$. We recall that faces
are the locally closed polyhedral cones (of all dimensions) into which $\Hc$ decomposes $\RR^n$.
They form a poset which we denote  $(\Cc=\Cc_\Hc,\leq)$. 

\vskip .2cm

 In this paper we propose an alternative description of $\Perv(\CC^n,\Hc)$ which is extremely
simple and appealing. It is given in terms of {\em matrix diagrams}, see Definition \ref{def:MD}, i.e.,
diagrams consisting of :

\begin{itemize}
\item[($0$)]  Vector spaces $E_{A,B}$ labelled by arbitrary pairs of faces $A,B\in\Cc$. 

\item[($1'$)]  A {\em representation} of $\Cc$  with respect
to the second argument, i.e., a transitive system of linear maps
\[
\del': E_{A,B_1}\lra E_{A,B_2}, \,\, B_1\leq B_2.
\]
\item[($1''$)]  An {\em anti-representation} of $\Cc$  with respect to the first argument, i.e.,  
 a transitive system of linear maps
\[
\del'': E_{A_2, B} \lra E_{A_1,B}, \,\,\, A_1\leq A_2.
\]
 \end{itemize}
 It is required that:
 \begin{itemize} 
 \item[($2$)] The maps $\del'$ and $\del''$  commute with each other, i.e., unite into a covariant functor
 $\Cc^\op\times\Cc\to\Vect_\k$.
 
 \item[($3$)]  If  the ``product cells''
 $iA+B_1$ and $iA+B_2$, $B_1\leq B_2$,
  lie in the same complex stratum, then the corresponding $\del'$
 is an isomorphism.
 Likewise, if $iA_1+B$ and $iA_2+B$, $A_1\leq A_2$, lie in the same complex stratum, 
 then the corresponding $\del''$ is an isomorphism. 
 \end{itemize}
 
 Using the ``Tits product'' 
 $ A\circ B$  on real cells (see \S \ref{sec:ex-comp}B below), one can give  
 (Remark \ref {rems:circle-cond}(a))
a reformulation of the condition (3) in terms involving only real cells.

\vskip .2cm
 
 Our main result, Theorem \ref {thm:main}, says that $\Perv(\CC^n,\Hc)$ is equivalent to the category
 of data $(E_{A,B},\del', \del'')$ satsfying the above conditions.  We also present, in Theorem 
 \ref{thm:affine}, a generalization to arrangements of affine hyperplane with real equations. 
 
 \vskip .2cm
 
 The simplest example for Theorem \ref {thm:main} is that of  $\Perv(\CC,0)$, the category of perverse sheaves on $\CC$
 with only singularity at $0$.  Our description identifies it with   the category of commutative
  $3\times 3$-diagrams below
 with the arrows at the outer  rim being isomorphisms: 
\be\label{eq:perv-C-0}
\xymatrix{
E_{-,+} \ar[r]^\simeq& E_{0,+}& \ar[l]_{\simeq} E_{+,+}
\\
E_{-,0}\ar[r]
\ar[u]^{\simeq} \ar[d]_{\simeq}
& E_{0,0}
\ar[u]\ar[d]
 & \ar[l] E_{+,0}
 \ar[u]_\simeq \ar[d]^\simeq
\\
E_{-,-}\ar[r]_\simeq  &E_{0,-} & \ar[l]^{\simeq} E_{+,-}
}
\ee

\paragraph*{}
Theorem  \ref {thm:main} is strikingly similar to the 
much more standard description of
 {\em constructible sheaves} on $(\CC^n, \Sc^{(0)})$ in terms of the quasi-regular cell decomposition
 of $\CC^n$ 
into the product cells $iA+B$. Such a sheaf $\Gc$ is given by its stalks $G_{A,B}$ at the $iA+B$
and generalization maps $\gamma', \gamma''$ just like in ($1'$) and ($1''$) but {\em covariant in both cases}. 
The condition that $\Gc$ is indeed $\Sc^{(0)}$-constructible means, just like in ($3$), that
$\gamma'$ and $\gamma''$ are isomorphisms whenever  the source and target correspond to product cells that lie
in the same complex stratum. 

Dually, an $\Sc^{(0)}$-{\em constructible 
cosheaf} (i.e., from the derived category point of view, a constructible complex Verdier dual to
a sheaf) is given by a diagram consisting of the $G_{AB}$ and maps  $\delta',\delta''$
{\em contravariant in both cases}, 
with the same properties. 

 Our result shows that perverse sheaves, occupying, intuitively, the middle position between sheaves and cosheaves,
admit a matching  description that is just as simple,  by reversing one of the two sets of arrows.

\vskip .2cm

Like a diagram describing a constructible sheaf, a matrix diagram has several ``layers''
(corresponding to complex strata $L_\CC^\circ$) with the property that the arrows within each layer are isomorphisms,
and therefore give a local system on the corresponding $L_\CC^\circ$ . For a constructiible sheaf $\Gc$
this local system is just the restriction of  $\Gc$ to ${L_\CC^\circ}$.  For a perverse sheaf 
$\Fc$ this corresponds to the restriction to $L_\CC^\circ$ of  the {\em hyperbolic restriction} of $\Fc$ to the closure of the stratum which is the  complex flat $L_\CC$,
see \cite{KS}, \S 5A. 
 The arrows between different layers describe the way such local systems are glued together. 
 
  For example, the outer rim of a diagram in \eqref{eq:perv-C-0}
represents a local system on $\CC\setminus\{0\}$ obtained by restricting 
the corresponding perverse sheaf from $\CC$
to $\CC\setminus\{0\}$, while the  full diagram can be seen as symbolically representing
the complex plane $\CC$ itself.  Further, the incoming and the outgoing arrows at the
middle term resemble the attractive and repulsive trajectories of a hyperbolic vector field
on $\CC=\RR^2$, very much in the spirit of the original philosophy of hyperbolic
localization \cite{GM-lefschetz, braden}.  

\paragraph*{}
Our method of proof of Theorem \ref{thm:main} is similar to that of \cite{KS}
but simplified, stripped, so to say, to the bare bones. As in \cite{KS}, the starting point is
the Cousin resolution $\Ec^\bullet(\Fc)$ of $\Fc\in\Perv(\CC^n,\Hc)$
associated to the system of tube cells $\RR^n+iA$, $A\in\Cc$, 
see \S \ref{sec:per-to-cous}. The matrix
diagram $(E_{A,B})$ corresponding to $\Fc$ is the linear algebra data
describing $\Ec^\bullet(\Fc)$ as a complex of cellular sheaves on the
cell decomposition   formed by the cells $iA+B$. That is,
the $E_{A,B}$ themselves are (up to sign factors) the stalks of the terms 
$\Ec^p=\Ec^p(\Fc)$.
The maps $\del'$ are the generalization maps describing the sheaf structure on each
$\Ec^p$. The maps $\del''$ describe the differentials $d: \Ec^p\to\Ec^{p+1}$. 
The  condition ($1'$) express the fact that each $\Ec^p$ is a sheaf. The condition
($1''$) expresses the requirement that $d^2=0$ in $\Ec^\bullet(\Fc)$. The remaining
conditions in ($2$) mean that $d$ is a morphism of cellular sheaves.
Thus any datum $(E_{A,B}, \del', \del'')$ satisfying ($1$)-($2$) 
always gives a  cellular complex $\Ec^\bullet$. 
The nontrivial part of Theorem \ref{thm:main} is that the condition ($3$) precisely
guarantees that this complex is in fact a perverse sheaf lying in $\Perv(\CC^n,\Hc)$,
in particular, that it is 
  constructible with respect to $\Sc^{(0)}$. 
  A crucial step here is a direct combinatorial identification of the Verdier dual to
   $\Ec^\bullet$
  (Proposition \ref{prop:E(F)-dual} and \S \ref{sec:vercousin}). The proof of 
  Theorem \ref{thm:main} is finished in \S \ref{sec:mat-to-per}.

\paragraph*{}

We would also  like to emphasize an important difference  between the present
description and that of 
\cite{KS}. While the approach 
  of    \cite{KS} is  centered around the
 ``real skeleton" $\RR^n\subset\CC^n$,
the linear algebra data in the present description are directly  tied to all the cells of a  cell decomposition  
of the underlying stratified space. Therefore the new approach can be viewed as somewhat bridging 
the gap between 
the {\em geometric definition} of perverse sheaves (via the t-structure on the derived
category) and various {\em combinatorial descriptions} (usually obtained by a
judicious choice of extra data).   Because of its ``local'' nature,
 it may be applicable in a wider range of situations
than just hyperplane arrangements.  A different  ``bridging'' approach, close to the ideas
of MacPherson \cite{macpherson:intersection}, is developed in \cite{DKSS}. 

\paragraph*{}
The research of M.K. was supported by World Premier International Research Center Initiative (WPI Initiative), 
 MEXT, Japan. 
 

\section{Generalities on arrangements}\label{sec:gen-arr}

We keep the notations and conventions of \cite{KS} which we recall for the reader's convenience.
Thus:

\vskip .2cm

 $\Hc$ is a finite arrangement of linear hyperplanes in $\RR^n$, assumed {\em central}, i.e.,
$\bigcap_{H\in\Hc} H = \{0\}$.  We choose, once and for all, a linear equation $f_H:\RR^n\to\RR$
for any $H\in\Hc$.  We denote 
\[
\sgn: \RR \lra \{+, -,0\}
\]
the standard sign function. 

\vskip .2cm

$(\Cc= \Cc_\Hc, \leq)$ is the poset of {\em faces} of $\Hc$, ordered by inclusion of the closures. 
By definition, $x,y\in\RR^n$ lie in the same face iff $\sgn  f_H(x)=\sgn  f_H(y)$
for each $H\in\Hc$. For a face $C$ and $H\in\Hc$ we denote 
\[
s_H(C)=\sgn\bigl(f_H|_C\bigr).
\]
The faces
are locally closed polyhedral subsets of $\RR^n$ forming a disjoint union (stratification) of $\RR^n$.
Each face is a topological cell, i.e., is homeomorphic to $\RR^p$ for some $p$. 

\vskip .2cm

If $A,B\in\Cc$, and $p>0$, the notation $A<_p B$ means that $A\leq B$ and $\dim(B) = \dim(A)+p$. 

\vskip .2cm

A {\em flat} of $\Hc$ is any intersection of hyperplanes from $\Hc$. This is understood to include
 $\RR^n$ itself
(the intersection of the empty set of hyperplanes). 

\vskip .2cm

For any subset $A\subset\RR^n$ we denote by $\LL(A)$ the $\RR$-linear span of $A$ and denote
\[
\pi_A: \RR^n\lra\RR^n/\LL(A)
\]
the projection.

\vskip .2cm

For any subspace $L\subset \RR^n$ we denote $L_\CC=L\otimes_\RR\CC\subset \CC^n$ its
complexification. In particular we have the complexified arrangement
$\Hc_\CC = \{H_\CC, H\in\Hc\}$ in $\CC^n$. 

\vskip .2cm

The space $\CC^n$ is equipped with the {\em complex stratification} ({\em $\CC$-stratification}
for short)  $\Sc^{(0)}$ whose strata are the generic parts of complexified flats, i.e., the
locally closed subvarieties ({\em $\CC$-strata})
\[
L_\CC^\circ \,=\, L_\CC\setminus \,\bigcup_{H\not\supset L} H_\CC,
\]
where $L$ is a flat of $\Hc$.  It also has the stratification (decomposition) $\Sc^{(2)}$ into the 
{\em product cells} $C+iD$, $C,D\in\Cc$. This decomposition refines $\Sc^{(0)}$. 

\vskip .2cm
We will also use the {\em intermediate} (or {\em Bj\"orner-Ziegler}) stratification  $\Sc^{(1)}$ of $\CC^n$ into strata $[C,D]$
parametrized by intervals $C\leq D$ in $\Cc$.  
 By definition (see   \cite{BZ} and \cite{KS} \S 2D),  $x+iy\in\CC^n$ lies in $[C,D]$, if:
 \be\label{eq:S1-cell}
 y\in C, \quad  
 x\in\pi_C^{-1}(\pi_C(D)).
\ee 
We consider each $\Sc^{(i)}$, $i=0,1,2$,  as a poset of strata ordered by inclusion of closures.

Ir is known (see \cite{KS} Prop. 2.10) that, denoting 
 by $\prec$ the relation of refinement of stratifications, we have
\[
\Sc^{(2)} \prec \Sc^{(1)} \prec \Sc^{(0)}. 
\]
In particular, we have the equivalence relations $\equiv_{\Sc^{(0)}}$ and $\equiv_{\Sc^{(1)}}$ on $\Sc^{(2)}=\Cc\times\Cc$
describing when the two product cells lie in one $\Sc^{(0}$-stratum (i.e., $\CC$-stratum) or in one $\Sc^{(1)}$-stratum. 
For simplicity we write $\equiv$ for $\equiv_{\Sc^{(0)}}$. We now describe these equivalence relations more explicitly,
starting with $\equiv$.

\vskip.2cm

For a face $C\in\Cc$ we denote 
\be\label{eq:H^C}
\Hc^C \,=\,\bigl\{H\in \Hc |\, H\supset C\}
\ee
the set of hyperplanes from $\Hc$ containing $H$. 

\begin{prop}\label{prop:A+iB=C+iD}
We have $A+iB \equiv C+iD$, if and only if 
\[
\Hc^A\cap\Hc^B \,=\, \Hc^C\cap\Hc^D.
\]
\end{prop}

\noindent{\sl Proof:} Let
\[
\Hc_\CC^{A+iB}\,=\,\bigl\{ H_\CC\in \Hc_\CC | \, H_\CC \supset A+iB\bigr\}. 
\]
 Then $A+iB \equiv C+iD$ if and only if $\Hc_\CC^{A+iB} = \Hc_\CC^{C+iD}$.  But for $H\in\Hc$, we have 
 $H_\CC\in \Hc_\CC^{A+iB}$ if and only if $H\in \Hc^A\cap \Hc^B$. 
Indeed, let $f=f_H: \RR^n\to\RR$ be the linear equaltion of $H$. Then
\[
f^\CC: \CC^n\lra \CC, \quad f^\CC(x+iy) = f(x) +i f(y), \,\,\, x,y\in\RR^n, 
\]
is a $\CC$-linear equaltion of $H_\CC$. So if $x,y\in\RR^n$ and $f^\CC(x+iy)=0$,
the $f(x)=f(y)=0$. \qed

\vskip .2cm

We now describe more explicitly the relation  $\equiv_{\Sc^{(1)}}$, i.e., the way how an $\Sc^{(1)}$-stratum
is decomposed into product cells. 

\begin{prop}\label{prop:S2vsS1}
Let $C\in\Cc$. Introduce an equivalence relation $\sim_C$ as the equivalence closure of the following relation
$\approx_C$: 
\[
B_1\approx_C B_2 \quad \text{iff } \quad B_1\leq B_2 \text{ and } B_1+iC \equiv B_2+iC. 
\]
Then, equivalence classes $\Bc$ of $\sim_C$ are on bijection with $\Sc^{(1)}$-strata $[C,D]$, $D\geq C$.
More precisely, each such stratum consists of product cells $B+iC$, $B\in\Bc$ for some $\sim_C$-class $\Bc$. 
\end{prop}

\noindent{\sl Proof:} By definition (the first condition in \eqref{eq:S1-cell}), each $[C,D]$ is the union of the
$B+iC$ where $B$ runs over some subset $\Bc\subset\Cc$. We prove that $\Bc$ is in fact an equivalence
class of $\sim_C$. 

We first prove that each such $\Bc$ is a union of equivalence classes of $\sim_C$. For this it suffices to
show that
\[
B_1\approx_C B_2 \quad \Rightarrow \quad B_1+iC \equiv_{\Sc^{(1)}} B_2+iC. 
\]
Indeed, suppose $B_1\leq B_2$ and $B_1+iC \equiv B_2+iC$. Then, in the notation of \eqref{eq:H^C},
\be\label{eq:B1-approx B_2}
\Hc^{B_1}\supset\Hc^{B_2}\quad\text{and} \quad \Hc^{B_1}\cap\Hc^C \,=\, \Hc^{B_2}\cap \Hc^C. 
\ee
The condition $B_1+iC \equiv_{\Sc^{(1)}} B_2+iC$ that we need to prove, means that $B_1$ and $B_2$
lie in the same region of the form $\pi_C^{-1}(\pi_C(D))$, i.e., that
\be\label{eq:sHB1=sHB2}
s_H(B_1) \,=\, s_H(B_2) \quad \text{for each } H\in\Hc^C.  
\ee
By \eqref {eq:B1-approx B_2}, for $H\in\Hc^C$ we have $s_H(B_1)=0$ iff $s_H(B_2)=0$. On the
other hand, since $B_1\leq B_2$, the difference between $s_H(B_1)$ and $s_H(B_2)$ can only be
that $s_H(B_2)\in \{\pm\}$ while $s_H(B_1)=0$. But this is impossible by  \eqref {eq:B1-approx B_2}. 
So  $B_1+iC \equiv_{\Sc^{(1)}} B_2+iC$ as claimed. 

\vskip .2cm

We now prove that $\Bc$ is a single equivalence class of $\sim_C$. For this we note that
the equivalence relation on $\Bc$ generated by $\leq$ (inclusion of closure of faces) has a single
equivalence class. Indeed,  $[C,D]$ is known to be
a cell (in particular,  connected) decomposed into product cells $B+iC$, $B\in\Bc$. So the relation
of inclusion of closures on these cells generates but a single equvalence class. 

But if 
 $B_1\leq B_2$ and $B_1+iC \equiv_{\Sc^{(1)}} B_2+iC$, then we have \eqref{eq:sHB1=sHB2}
and so \eqref{eq:B1-approx B_2} so $B_1\approx_C B_2$.  \qed

 
\section{Generalities on cellular sheaves and perverse sheaves}

We fix a base field $\k$ and denote $\Vect_\k$ the category of $\k$-vector spaces.
By a {\em sheaf} we always mean a sheaf of $\k$-vector spaces. 
For a topological space $X$ we denote by $\Sh_X$ the category of
sheaves  on $X$ and by $D^b\Sh_X$ the bounded derived category
of $\Sh_X$.  For $V\in\Vect_\k$ we denote by $\ul V_X$ the constant sheaf on $X$ with
stalk $V$. 

\vskip .2cm

If $(X,\Sc)$ is a stratified space, then we denote $\Sh_{X,\Sc}$ the category of
$\Sc$-{\em constructible} sheaves on $X$, i.e., sheaves $\Fc$ which are locally constant,
with finite-dimensional stalks, on each stratum of $\Sc$. We denote by $D^b_\Sc \Sh_X\subset
D^b\Sh_X$ the full subcategory of {\em   $\Sc$-constructible complexes}, i.e.,
of complexes $\Fc$ such that each cohomology sheaf $\ul H^i(\Fc)$ lies in
$\Sh_{X,\Sc}$. The triangulated category $D^b_\Sc \Sh_X$ has a perfect dualty,
the {\em Verdier duality}, denoted $\Fc\mapsto\DD(\Fc)$. 

\vskip .2cm

By a {\em cell} we mean a topological space $\sigma$ homeomorphic to $\RR^d$ for some $d$.
A {\em cellular space} is a stratified space $(X.\Sc)$ such that each stratum is a cell.
We consider $\Sc$ as the poset of cells with the order $\leq$ given bu inclusion of the closures. 
A {\em cellular sheaf}, resp. {\em cellular complex} on a cellular space $(X,\Sc)$ is an
$\Sc$-constructible sheaf, resp. complex on $X$. For such a sheaf, resp. complex $\Fc$
and a cell $j_\sigma:\sigma\hra X$ we denote by
\[
\Fc\para|_\sigma \,=\, R\Gamma(\sigma, j_\sigma^*\Fc)
\]
the stalk of $\Fc$ at $\sigma$. 
We recall from \cite{KS} \S 1D the concept of a {\em quasi-regular cellular space}
as well as the following fact. 

\begin{prop}\label{prop:cell-sheaves}
Let $(X, \Sc)$ be a quasi-regular cellular space. Then $\Sh(X,\Sc)$ is equivalent to the category of
representations of the poset $(\Sc,\leq)$, i.e., of data formed by:
\begin{itemize}
\item[(0)] Vector spaces $G_\sigma, \sigma\in \Sc$.

\item[(1)] Linear maps $\gamma_{\sigma,\sigma'}: G_\sigma\to G_{\sigma'}$ given for each
$\sigma\leq\sigma'$ and satisfying the transitivity relations:
\item[(2)] $\gamma_{\sigma, \sigma''} = \gamma_{\sigma', \sigma''}\circ \gamma_{\sigma, \sigma'}$ for
any $\sigma\leq\sigma'\leq\sigma''$. 
\end{itemize}
Explicitly, to a sheaf $\Gc\in\Sh(X,\Sc)$ there corresponds the datum formed by the stalks
$G_\sigma=\Gc\para_{\sigma}$ and the generalization maps
$\gamma_{\sigma,\sigma'}:\Gc\para_\sigma\to\Gc\para_{\sigma'}$. \qed
\end{prop}

We specialise to the situation of \S\ref{sec:gen-arr} and take  $X=\CC^n$. The stratifications $\Sc^{(1)}$ and $\Sc^{(2)}$
are quasi-regular cell decompositions of $\CC^n$, while $\Sc^{(0)}$ is not.

We will use the involution
\be\label{eq:tau}
\tau:\CC^n\to\CC^n, \quad (x+iy) \,  \mapsto\,  (y+ix). 
\ee
This involution is not $\CC$-linear but it preserves $\Hc$ and the stratifications $\Sc^{(0)}$
(stratum by stratum, i.e., $\tau$ preserves each stratum)  and $\Sc^{(2)}$ (as a whole, i.e., $\tau$ takes each stratum to another stratum).   But it does  not preserve
$\Sc^{(1)}$. We denote by $\tau\Sc^{(1)}$ the new stratification of $\CC^n$ whose strata are
obtained by applying $\tau$ to the strata
of $\Sc^{(1)}$. 

\begin{prop}\label{propS1+tS1}
Let $\Fc\in\Sh(\CC^n, \Sc^{(2)})$. Suppose that $\Fc$ is both $\Sc^{(1)}$-constructible and $\tau\Sc^{(1)}$-constructible.
Then $\Fc$ is $\Sc^{(0)}$-constructible. 
\end{prop}

\noindent{\sl Proof:} By definition, being $\Sc^{(0)}$-construcible means that for inclusion $B_1+iC_1 \leq B_2+iC_2$
of (closures of) product cells such that $B_1+iC_1 \equiv B_2+iC_2$, the corresponding generalization map
\[
\gamma_{B_1+iC_1, B_2+iC_2}: \Fc\para_{B_1+iC_1} \lra  \Fc\para_{B_2+iC_2}
\]
 is an isomorphism. Now, $\Sc^{(2)} = \Cc\times\Cc$ is the product stratification, and, moreover,
 the stratification of each  complexified flat of $\Hc$ into $\Sc^{(2)}$-strata is also a product
 stratification. So it suffices to prove the isomorphicity of $\gamma_{B_1+iC_1, B_2+iC_2}$
 in two separate cases, horizontal and vertical:
 \begin{itemize}
\item [(1)] $B_1\leq B_2$, $C_1=C_2=C$, and  $B_1+iC \equiv B_2+iC$. 

\item[(2)] $B_1=B_2=B$, $C_1\leq C_2$ and $B+iC_1\equiv B+iC_2$. 
 \end{itemize}
Now,  inclusions of type (1) generate, by Proposition \ref{prop:S2vsS1}, inclusions of $\Sc^{(2)}$-cells
within the same $\Sc^{(1)}$-cell. Similarly, inclusions of type (2) generate inclusions of $\Sc^{(2)}$-cells
within the same $\tau \Sc^{(1)}$-cell. So if  $\Fc$ is both $\Sc^{(1)}$-constructible and $\tau\Sc^{(1)}$-constructible,
then the generalization maps corresponding to (1) and (2) are all isomorphisms hence $\Fc$ is
$\Sc^{(0)}$-constructible. \qed

\vskip .2cm

We denote by $\Perv(\CC^n,\Hc)\subset D^b_{\Sc^{(0)}} \Sh_{\CC^n}$   the category of perverse sheaves
(with respect to the middle perversity),
 which are are $\Sc^{(0)}$-constructible. 
Explicitly, we normalize the perversity conditions by saying that $\Fc$ is perverse, if:
\begin{itemize}
\item[($P^-$)] The sheaf $\ul H^p(\Fc)$ is supported on a subvariety of complex codimension
$\geq p$. 

\item [($P^+$)] If $S$ is a stratum of $\Sc^{(0)}$ of complex codimension $p$, then the
sheaves $\ul \HH^i_S(\Fc)$ of hypercohomology with support in $S$, are zero for $i<p$. 
\end{itemize}
In this normalization, a perverse sheaf reduces, on the open stratum, to a local system
in degree $0$. 
As well known , the condition ($P^+$) for $\Fc$ is equivalent to ($P^-$) for the shifted
Verdier dual
\[
\Fc^* \,=\, \DD(\Fc) [-2n]. 
\]
 The functor $\Fc\mapsto \Fc^*$ is thus a 
   perfect duality  on $\Perv(\CC^n,\Hc)$. Since the involution $\tau$ preserves $\Sc^{(0)}$,
   the pullback functor $\tau^{-1}$ preserves the category $\Perv(\CC^n,\Hc)$ and so
   this category has another perfect duality
   \be\label{eq:F-tau}
   \Fc\,\mapsto \, \Fc^\tau := \tau^{-1} \Fc^*. 
   \ee
   
   
   \section{Matrix diagrams and the main result}

   \begin{Defi}\label{def:MD}
   A {\em matrix diagram} of type $\Hc$ is a collection $E$ of the following data:
   \begin{itemize}
   \item[($M0$)] Finite-dimensional $\k$-vector spaces $E_{A,B}$ given for any two faces $A,B\in\Cc$.
   
   \item[($M1$)] Linear maps
   \[
   \del'=\del'_{(A|B_1, B_2)}: E_{A,B_1} \lra E_{A, B_2},
   \text{ given for any faces $A$ and $B_1\leq B_2$,  }
   \]
      \[
   \del'' = \del''_{(A_2, A_1|B)}: E_{A_2, B}\lra E_{A_1, B}, 
   \text{ given for any faces $A_1\leq A_2$ and $B$,}
   \]
  satisfying the following conditions:
   
   \item[($M2$)] The maps $\del', \del''$ define a representation of the poset $\Cc^\op\times\Cc$ in
   $\Vect_\k$. That is, we have
   \[
   \begin{gathered}
   \del'_{(A|B_2, B_3)}\circ \del'_{(A|B_1, B_2)} \,=\, \del'_{(A|B_1, B_3)}, \quad \forall \,\, A 
   \text{ and } B_1\leq B_2\leq B_3;
   \\
   \del''_{(A_2, A_1|B)}\circ\del''_{(A_3, A_2|B)} \,=\,
   \del''_{(A_3, A_1|B)}, \quad \forall \,\, A_1\leq A_2\leq A_3 \text{ and } B; 
   \\
   \del'_{(A_2, A_1|B_2)} \circ \del''_{(A_2|B_1, B_2)} \,=\, \del''_{(A_1|B_1, B_2)}\circ
   \del'_{(A_2, A_1|B)},\quad\forall \,\, A_1\leq A_2 \text{ and } B_1\leq B_2. 
   \end{gathered}
   \]
   
   \item[($M3'$)] If $A$ and $B_1\leq B_2$ are such  $B_1+iA\equiv B_2+iA$ (that is, these product cells lie
   in the same complex stratum),
   then $\del'_{(A|B_1, B_2)}$ is an isomorphism.
   
   \item[($M3''$)] Similarly, if $A_1\leq A_2$ and $B$ are such that $B+iA_1\equiv B+iA_2$,
    then $\del''_{(A_2, A_1|B)}$ is an isomorphism. 
   
   \end{itemize}
   \end{Defi}
   
   We denote by $\Mc_\Hc$ the  category of matrix diagrams of type $\Hc$. This category
   is abelian and has a perfect duality
   \be\label{eq:dual-matrix}
   E\mapsto E^*, \quad (E^*)_{A,B} = (E_{B,A})^*
   \ee
(``hermitian conjugation'') with the maps $\del'$ for $E^*$ being the dual of the $\del''$ for $E$
and the $\del''$ for $E^*$ being dual to the $\del'$ for $E$.

\begin{rem}\label{rem:A+ib=B+iA}
We  note that any matrix diagram is ``symmetric'' in the following weak sense. Since
$B+iA$ and $A+iB$ lie in the same complex stratum by Proposition 
\ref{prop:A+iB=C+iD}, 
 we have  an isomorphism (non-canonical) $E_{A,B}\simeq E_{B,A}$. It is given by the monodromy
of the ``layer'' (system of isomorphic maps $\del', \del''$) 
of the matrix diagram  corresponding to this complex stratum. 
\end{rem}

Our main result is as follows.

\begin{thm}\label{thm:main}
We have mutually quasi-inverse equivalences of categories
 \[
  \xymatrix{
  \Perv(\CC^n, \Hc)  \ar@<.5ex>[r]^{\hskip 0.7cm \EE}& \Mc_\Hc \ar@<.5ex>[l]^
 {\hskip .7cm \GG  } 
  }
  \]
  taking the  twisted Verdier duality $\Fc\mapsto\Fc^\tau$ to the duality \eqref {eq:dual-matrix}. 
\end{thm}


\section{From a perverse sheaf to a matrix diagram:
 the Cousin complex}\label{sec:per-to-cous}

Here we construct a functor $\EE: \Perv(\CC^n,\Hc)\to\Mc_\Hc$. We use the 
general analysis of perverse sheaves on $(\CC^n,\Hc)$ in terms of their Cousin complexes
\cite{KS}.

For a face $A\in\Cc$ we denote $\lambda_A: \RR^n+iA \hookrightarrow\CC^n$ 
the embedding of the corresponding ``tube cell".

\begin{prop}\label{prop:cousin-properties}
Let $\Fc\in\Perv(\CC^n,\Hc)$. Then:

\vskip .2cm

(a) The complex $\lambda_A^!\Fc$ reduces to a single sheaf $\wt\Ec_A = \wt\Ec_A(\Fc)$ 
in degree $\cod(A)$.

\vskip .2cm

(b) The complex $\lambda_{A*} \wt\Ec_A(\Fc)$ also reduces to a single sheaf $\Ec_A=\Ec_A(\Fc)$
supported on $\RR^n+i\ol A$. 

\vskip .2cm

(c) $\Ec_A(\Fc)$, considered as a sheaf on $\RR^n+i\ol A$, is the pullback, with respect to the
projection to $\RR^n$, of a sheaf on $\RR^n$, constructible with respect to the stratification $\Cc$. 
\end{prop}

\noindent{\sl Proof:} For $A=\{0\}$, this is Prop. 4.9(a) of \cite{KS}. For an arbitrary $A$ this
follows by further applying Prop. 3.10 and  Cor. 3.22 from \cite{KS}. \qed 

\vskip .2cm
Further, the standard coboundary maps on the sheaves of cohomology with support
give the {\em Cousin complex}
\be
\Ec^\bullet(\Fc) \,=\,\biggl\{ \bigoplus_{\cod(A)=0} \Ec_A(\Fc)\buildrel\delta\over
\to  \bigoplus_{\cod(A)=1} \Ec_A(\Fc)\buildrel\delta\over\to\cdots \buildrel\delta\over\to
\Ec_0(\Fc)\biggr\}
\ee
which is a complex of sheaves on $\CC^n$ canonically isomorphic to $\Fc$ in $D^b\Sh_{\CC^n}$,
see \cite{KS}, Cor. 4.11. The grading in $\Ec^\bullet(\Fc) $ is by $\cod(A)$. 

Further, the matrix elements of $\delta$ which are morphisms of sheaves
\[
\delta_{A_2, A_1}: \Ec_{A_2}(\Fc) \lra\Ec_{A_1}(\Fc), \quad \cod(A_2) = \cod(A_1)-1, 
\]
are nonzero only if $A_1<_1 A_2$. The condition $\delta^2=0$ means that the
$\delta_{A_2, A_1}$ anticommute with each other. That is, for any faces
$A_1, A_2\neq A'_2, A_3$ such that $A_1\leq_1 A_2, A'_2\leq_1 A_3$
(a commutative square in $\Cc$ as a category), we have
\[
\delta_{A_2, A_1}\circ\delta_{A_3, A_2} \,=\,- \delta_{A'_2, A_1}\circ\delta_{A_3, A'_2}:
\Ec_{A_3}(\Fc)\lra \Ec_{A_1}(\Fc). 
\]
This anticommutativity can be converted to commutativity in a standard way by
``introducing signs''. More precisely, for any cell   $\sigma$
 let
\[
\orr(\sigma ) \,=\, H^{\dim(\sigma)}_c(\sigma,\k)
\]
be the $1$-dimensional orientation vector space of $\sigma$.  Note that $\orr(\sigma)^{\otimes 2}$
is canonically identified with $\k$. In particular, every face $A$ being a cell, we have the
space $\orr(A)$. 
  For any $A_1 <_1 A_2$ we have a canonical isomorphism
  \[
  \psi_{A_1, A_2}: \orr(A_1)\lra \orr(A_2)
  \]
which is the matrix element of the differential in the cellular cochain complex of $\ol A_2$
with coefficients in $\k$.  For  $A_1\leq_1 A_2, A'_2\leq_1 A_3$ as above, the isomorphisms
$\psi$ anticommute. So we get the following:

\begin{prop}\label{prop:del-commute}
(a) The morphisms
\[
\del_{A_2, A_1} =\delta_{A_2, A_1}\otimes\psi_{A_1, A_2}^{-1}: \Ec_{A_2}(\Fc)\otimes_\k\orr(A_2)
\lra\Ec_{A_1}(\Fc)\otimes_\k \orr(A_1)
\]
satisfy the commutativity constraints. That is, for any $A_1\leq_1 A_2, A'_2\leq_1 A_3$ as above,
\[
\del_{A_2, A_1}\circ\del_{A_3, A_2} \,=\,\del_{A'_2, A_1}\circ\del_{A_3, A'_2}. 
\]
(b) The maps $\del_{A_2, A_1}$, $A_1<_1 A_2$ extend to a representation of the poset
$\Cc^\op$ in $\Sh_{\CC^n}$ which takes $A$ to $\Ec_A(\Fc)\otimes \orr(A)$. 
In other words, for any $A_1<_p A_2$, $p\geq 1$,  we have a morphism of sheaves
\[
\del_{A_2, A_1}: \Ec_{A_2}(\Fc)\otimes \orr(A_2)
\lra\Ec_{A_1}(\Fc)\otimes \orr(A_1)
\]
defined as
\[
\del_{A_2, A_1} \,=\,\del_{A'_1, A_1}\circ \del_{A'_2, A'_1}\circ\cdots\circ \del_{A_2, A'_{p-1}},
\]
for  any chain $A_1 <_1 A'_1 <_1 \cdots <_1 A'_{p-1} <_1 A_2$, the result being independent
on the choice of such chain. These morphisms satisfy the transitivity condition for any
$A_1 <_p A_2 <_q A_3$. 
\qed
\end{prop}

\vskip .2cm

Let now $A,B\in\Cc$ be two faces. We associate to $\Fc\in\Perv(\CC^n,\Hc)$ the vector space
\[
E_{A,B}= E_{A,B}(\Fc)\, := \, (\Ec_A(\Fc)\otimes\orr(A))_{B+i0}, 
\]
 the stalk of $\Ec_A(\Fc)\otimes\orr(A)$ at the cell
$B+i0\subset \RR^n+iA$. Because of Proposition \ref{prop:cousin-properties}(c), we have
a canonical identification
\be\label{eq:EAB-indentif}
E_{A,B}(\Fc) \,\simeq \,  (\Ec_A(\Fc)\otimes\orr(A))_{B+iA'}, \quad A'\leq A
\ee
with the stalk at  $B+iA'$ for any $A'\leq A$. 

\vskip .2cm

If we have faces $A$ and $B_1\leq B_2$, then we define
\[
\del'_{(A|B_1, B_2)}: E_{A, B_1}(\Fc)\lra E_{A, B_2}(\Fc)
\]
to be the generalization map of the cellular sheaf $\Ec_A(\Fc)\otimes\orr(A)$ from $B_1+i0$
to $B_2+i0$.

If we have faces $A_1\leq A_2$ and $B$, then we define
\[
\del''_{(A_1, A_2|B)}: E_{A_2,B}(\Fc)\lra E_{A_1,B}(\Fc)
\]
to be the map of stalks at $B+i0$ induced by the morphism of sheaves $\del_{A_2, A_1}$. 

\begin{prop}\label{prop:E=MD}
  The system $E(\Fc)=(E_{A,B}(\Fc), \del', \del'')$
is a matrix diagram of type $\Hc$. We have therefore a functor
\[
\EE: \Perv(\CC^n,\Hc)\lra \Mc_\Hc,\quad \Fc\mapsto E(\Fc). 
\]
\end{prop}

\noindent{\sl Proof:} We first prove the conditions ($M2$) of Definition \ref{def:MD} of a matrix diagram.
The first condition in ($M2$) follows from the fact that $\Ec_A(\Fc)\otimes\orr(A)$ is a cellular
sheaf amd so its generalization maps are transitive. The second condition in ($M2$) follows
from Proposition \ref{prop:del-commute}(b). Finally, the third condition in ($M2$) follows
from the fact that $\del_{A_2, A_1}$ is a morphism of cellular sheaves and so the maps
it induces on the stalks, commute with the generalization maps. 

Let us prove the condition ($M3'$) of Definition \ref{def:MD}. By construction, $\Ec_A=\Ec_A(\Fc)$
is locally constant on the intersection of each stratum of $\RR^n+i\ol A$ (i.e., of each
$\RR^n+iA'$, $A'\leq A$ with each $\CC$-stratum. So if $B_1+iA \equiv B_2+iA$
  and $B_1\leq B_2$, then the generalization map on the
stalks
\[
\gamma_{B_1+iA, B_2+iA}: (\Ec_A)_{B_1+iA} \lra (\Ec_A)_{B_2+iA}
\]
is an isomorphism. But in virtue of \eqref{eq:EAB-indentif}, this map is identified with
\[
\del'_{(A|B_1, B_2)}: E_{A, B_1} \lra E_{A, B_2}
\]
and so the latter map is an isomorphism, proving ($3'$). 

\vskip .2cm

The property ($M3''$) for $E(\Fc)$ will follow from ($M3'$) for $E(\Fc^\tau)$ if we prove the
 following fact.

\begin {prop}\label{prop:E(F)-dual}
The system $E(\Fc^\tau)$ is identified with the dual system to $E(\Fc)$, that is,
$E_{A,B}(\Fc^\tau)$ is  identified with  $(E_{B,A}(\Fc))^*$ so that the maps $\del'$ (resp. $\del''$)
for $E(\Fc^\tau)$ are identified with the duals  to the $\del''$ (resp. $\del'$) for $E(\Fc)$. 
\end{prop}

This will be done in the next section. 


\section{Verdier duality and the Cousin complex}\label{sec:vercousin}

In this section we rewrite the Cousin complex $\Ec^\bullet(\Fc)$ in a way manifestly
compatible with Verdier duality. We start by  one more general remark on cellular sheaves.
 .

\vskip .2cm

Let $(X,\Sc)$ be a quasi-regular cellular space with cell embeddings denoted $j_\sigma:\sigma\hra X$. 
  Let $\Gc$ be a cellular sheaf
 on $X$ given by the   linear algebra data $(G_\sigma, \gamma_{\sigma,\sigma'})$ of Proposition \ref{prop:cell-sheaves}. 
 Then, these data give a complex in the derived category $D^b\Sh_X$:
\be\label{eq:cell-resol-der}
\bigoplus_{\dim(\sigma)=0} j_{\sigma !} \ul{G_\sigma}_\sigma \lra 
\bigoplus_{\dim(\sigma)=1} j_{\sigma !} \ul{G_\sigma}_\sigma[1]  \lra\cdots
\ee
whose total object is $\Gc$, see \cite{KS} (1.12). We say that \eqref{eq:cell-resol-der}
is a resolution of $\Gc$. Note that given just vector spaces $G_\sigma$, the datum of
such a complex is equivalent to the datum of transitive $\gamma_{\sigma,\sigma'}$,
i.e., of a cellular sheaf with these stalks.

\vskip .2cm

We apply this to the sheaf $\Ec_A(\Fc)$ on the cellular space formed by $\CC^n$ with
the stratification $\Sc^{(2)}$ into product cells $B+iA$. Given any such cell, we have
a commutative diagram of embeddings
\[
\xymatrix{
B+iA 
\ar[dr]^{\eps^{BA}}
\ar[d]_{l^{BA}}
\ar[r]^{k^{BA}} & \RR^n+iA
\ar[d]^{\lambda_A}
\\
B+i\RR^n \ar[r]_{\kappa_B}& \CC^n. 
}
\]

\begin{prop}
let $V$ be a $\k$-vector space. Then we have a canonical identification
\[
\lambda_{A*}\,  k^{BA}_! \, \ul V_{B+iA} \,\simeq \, \kappa_{B!}\,  l^{BA}_* \, \ul V_{B+iA}. 
\]
\end{prop} 
\noindent {\sl Proof:} The stalk of either of these sheaves at $x+iy\in\CC^n$ is
\[
\begin{cases}
V, &\text{ if } x\in B \text{ and } y\in\ol A,
\\
0,& \text{ otherwise}. \hskip 5cm \qed
\end{cases}
\]

\vskip .2cm

We will denote the sheaf in the proposition by $\eps^{BA}_{!*}\, \ul V_{B+iA}$
and refer to it as a {\em 1-cell sheaf of $(!*)$-type}. We similarly denote 1-cell sheaves
of $(*!)$-type as
\[
\eps^{BA}_{*!} \, \ul V_{B+iA} \,=\, 
\lambda_{A!}\,  k^{BA}_* \, \ul V_{B+iA} \,\simeq \, \kappa_{B*}\,  l^{BA}_! \, \ul V_{B+iA}. 
\]

\begin{prop}
Let $\Fc\in\Perv(\CC^n,\Hc)$ and $A\in\Cc$. Then the sheaf $\Ec_A(\Fc)$ has a resolution
(in $D^b\Sh_{\CC^n}$) of the form
\[
\bigoplus_{\dim(B)=0} \eps^{BA}_{!*}\ul {E_{AB}}_{B+iA} \lra 
\bigoplus_{\dim(B)=1} \eps^{BA}_{!*}\ul {E_{AB}}_{B+iA}[1] \lra \cdots
\]
with the differentials given by the maps $\del'$. 
\end{prop}

\noindent {\sl Proof:} This is an instance of \eqref{eq:cell-resol-der}. 
 It simply reflects the fact that $\Ec_A$ is the sheaf on $\CC^n$ coming from the sheaf
 on $\RR^n+i\ol A$ which is pulled back from the $\Cc$-constrructible (cellular) sheaf on
 $\RR^n$ with stalks $E_{AB}$ and generalization maps $\del'$. \qed
 
 \begin{cor}\label{cor:E(F)-bicell}
 Any $\Fc\in\Perv(\CC^n,\Hc)$ has a canonical resolution (in $D^b\Sh_{\CC^n}$) in the form
 of the double complex
 \[
 \xymatrix{
 \bigoplus\limits_{\cod(A)=0\atop \dim(B)=0}  \eps^{BA}_{!*}\, \ul{E_{AB}\otimes\orr(A)}_{B+iA} 
 \ar@<-10.5ex>[d]
 \ar[r]&  \bigoplus\limits_{\cod(A)=0\atop \dim(B)=1}  \eps^{BA}_{!*}\, \ul{E_{AB}\otimes\orr(A)}_{B+iA} [1]
 \ar@<-11.5ex>[d]
 \ar[r]&\cdots
 \\
 \bigoplus\limits_{\cod(A)=1\atop \dim(B)=0}  \eps^{BA}_{!*}\, \ul{E_{AB}\otimes\orr(A)}_{B+iA} 
 \ar[r]
 \ar@<-10.5ex>[d]
 &  \bigoplus\limits_{\cod(A)=1\atop \dim(B)=1}  \eps^{BA}_{!*}\, \ul{E_{AB}\otimes\orr(A)}_{B+iA} [1]
 \ar[r]
 \ar@<-11.5ex>[d]&\cdots
 \\
 &&
 \\
\hskip -3.8cm  \vdots & \hskip -3.9cm \vdots &
 }
 \]
with the horizontal differentials given by the maps $\del'$ and the vertical differentials given
by the $\del''$. 
 \end{cor}
 
 \noindent{\sl Proof:} This is just the Cousin complex written in terms of $1$-cell sheaves
 (of $(!*)$-type). \qed
 
 \vskip .2cm
 
 We now prove Proposition \ref{prop:E(F)-dual}. For this, we apply the shifted Verdier duality
 to the double complex in Corollary \ref{cor:E(F)-bicell} and note the three standard facts:
 \begin{itemize}
 \item $\DD$ interchanges $f_!$ and $f_*$.
 
 \item  For a cell $\sigma$ of real dimension $d$ and a finite-dimensional $\k$-vector
 space $V$,  we have $\DD(\ul V_\sigma)= \ul{V^*\otimes \orr(\sigma)}_\sigma [d]$. 
 
 \item $\orr(\sigma)^{\otimes 2}\simeq \k$ canonically. 
 \end{itemize}
 We conclude that $\Fc^*$ has a resolution in $D^b\Sh_{\CC^n}$ in the form of the double
 complex
 \be\label{double-resol-dual}
 \xymatrix{
 \bigoplus\limits_{\dim(A)=0\atop \cod(B)=0}  \eps^{BA}_{*!} \, \ul{E^*_{AB}\otimes\orr(B)}_{B+iA} 
 \ar@<-10.5ex>[d]
 \ar[r]&  \bigoplus\limits_{\dim(A)=0\atop \cod(B)=1}  \eps^{BA}_{*!} \, \ul{E^*_{AB}\otimes\orr(B)}_{B+iA}  
 \ar@<-11.5ex>[d]
 \ar[r]&\cdots
 \\
 \bigoplus\limits_{\dim(A)=1\atop \cod(B)=0}  \eps^{BA}_{*!}\, \ul{E^*_{AB}\otimes\orr(B)}_{B+iA} [1]
 \ar[r]
 \ar@<-10.5ex>[d]
 &  \bigoplus\limits_{\cod(A)=1\atop \dim(B)=1}  \eps^{BA}_{!*}\, \ul{E_{AB}\otimes\orr(A)}_{B+iA} [1]
 \ar[r]
 \ar@<-11.5ex>[d]&\cdots
 \\
 &&
 \\
\hskip -3.8cm  \vdots & \hskip -3.9cm \vdots &
 }
\ee
 with the horizontal differentials given by the duals to the $\del''$ for $E(\Fc)$ and the vertical
 differentials given by the duals of the $\del'$ for $E(\Fc)$. It corresponds, therefore,
 to the dual system $E(\Fc)^*$. 
 
 On the other hand, we can form the Cousin resolution of $\Fc^*$ but using the real, not imaginary
 tube cells $\kappa_B: B+i\RR^n\hookrightarrow\CC^n$. This gives the sheaves 
 \[
 \wc\Ec_B(\Fc^*) \,=\,\kappa_{B*}\,\kappa_B^!\, \Fc^*[\cod (B)]
 \]
 and the resolution
 \[
 \wc\Ec^\bullet(\Fc^*) \,=\,\biggl\{ \bigoplus_{\cod(B)=0} \wc\Ec_B(\Fc^*) \lra
 \bigoplus_{\cod(B)=1} \wc\Ec_B(\Fc^*) \lra
\cdots\biggr\}
 \]
of $\Fc^*$. Writing out each $\Ec_B(\Fc^*)$ in terms of $1$-cell sheaves of type $(*!)$,
we get a double complex of the form  \eqref{double-resol-dual} which is a resolution of
$\Fc^*$. But the Cousin resolution of $\Fc^*$ with respect to the cells $B+i\RR^n$ is the same
as the Cousin resolution of $\tau^{-1}\Fc^*$ with respect to the the cells $\RR^n+iA$. We conclude
that the complex \eqref{double-resol-dual}, associated to $E(\Fc)^*$
must reduce, after applying $\tau$, to the complex of  Corollary \ref{cor:E(F)-bicell}
 describing $E(\Fc^*)$. This means that the linear algebra data underlying the two
 complexes must be identified, i.e., $E(\Fc)^*\simeq E(\Fc^\tau)$. 
 
 This finishes the proof of Propositions   \ref{prop:E(F)-dual} and \ref{prop:E=MD}. 
 

\section{From a matrix diagram to a perverse sheaf}\label{sec:mat-to-per}

We now construct a functor
\[
\GG: \Mc_\Hc\lra\Perv(\CC^n, \Hc)
\]
by reversing the procedure used to extract the matrix diagram $E(\Fc)$ from the Cousin
complex $\Ec^\bullet(\Fc)$. 

Let $E=(E_{A,B}, \del', \del'')\in\Mc_\Hc$ be given. For each face $A\in\Cc$ we form the
cellular shef $\Ec_A=\Ec_A(E)$ on $(\CC^n, \Sc^{(2)})$ which is supported on $\RR^n+i\ol A$
and pulled there from the cellular sheaf on $(\RR^n,\Cc)$ with stalks $E_{A,B}\otimes\orr(A)$
and generalization maps $\del'\otimes \Id$. In other words, $\Ec_A$ is
constant on each $B+i\ol A$ and
\be\label{eq:EA on closed}
\Ec_A|_{B+i\ol A} \,=\,\ul{ E_{A,B}\otimes\orr(A)}_{B+i\ol A} 
\ee
Further, the commuting maps $\del''$ in $E$ give, after tensoring with the $\orr(A)$,
anticommuting morphisms of sheaves
\[
\delta_{A_2, A_1}: \Ec_{A_2}(E) \lra\Ec_{A_1}(E),\quad A_1 <_1 A_2,
\]
and so we can form the complex of sheaves
\[
\Ec^\bullet(E) \,=\,\biggl\{ \bigoplus_{\cod(A)=0} \Ec_A(E) \buildrel \delta\over\lra 
 \bigoplus_{\cod(A)=1} \Ec_A(E) \buildrel \delta\over\lra \cdots\biggr\}. 
\]

\begin{prop}\label{prop:E(E)S0constr}
The complex $\Ec^\bullet(E)$ is  $\Sc^{(0)}$-constructible.
\end{prop}

\noindent{\sl Proof:} By definition, $\Ec^\bullet=\Ec^\bullet(E)$ is $\Sc^{(2)}$-constructible. By
Proposition
\ref{propS1+tS1}, it suffices to prove that it is both $\Sc^{1)}$-constructible and $\tau\Sc^{(1)}$-constructible.

Let us first prove that $\Ec^\bullet$ is $\Sc^{1)}$-constructible. 
By Proposition \ref{prop:S2vsS1} this is equivalent to the following condition:

\begin{itemize}
\item[(Q)] If $C_1\leq C_2$ and $D$ are such that $C_1+iD\equiv C_2+iD$, then the generalization map
\[
\gamma_{C_1+iD, C_2+iD}: \Ec^\bullet\para_{C_1+iD} \lra  \Ec^\bullet\para_{C_2+iD} 
\]
is a quasi-isomorphism of complexes. 
\end{itemize}

 We claim that 
  $\gamma_{C_1+iD, C_2+iD}$ is in fact an isomorphism, not just a quasi-isomorphism
of complexes. More precisely, we claim that for any summand $\Ec_A$ of $\Ec^\bullet$,
the corresponding generalization map
\[
\gamma^{\Ec_A}_{C_1+iD, C_2+iD}: \Ec_A\para_{C_1+iD} \lra  \Ec_A\para_{C_2+iD} 
\]
is an isomorphism of vector spaces. To see this, note that by construction, see
\eqref{eq:EA on closed}, we have for any $C,D$:
\be
\Ec_A\para_{C+iD} \,=\,\begin{cases}
E_{A,C}\otimes\orr(A),& \text{if } D\leq A,
\\
0,&\text{ otherwise.} 
\end{cases}
\ee
So
\[
\gamma^{\Ec_A}_{C_1+iD, C_2+iD}\,=\,\begin{cases}
\del'\otimes\Id: E_{A, C_1}\otimes\orr(A)\to E_{A,C_2}\otimes\orr(A),& \text{ if } D\leq D, 
\\
0: 0\lra 0,& \text{ otherwise}. 
\end{cases}
\]
So the fact that it is an isomorphism, follows from condition (3') of Definition 
\ref{def:MD} of matrix diagram and the next lemma.

\begin{lem}\label{lem:C-stratum-enlarge}
Let $C_1, C_2$ and $D$ be such that $C_1+iD\equiv C_2+iD$ lie in the same complex
stratum. Let $A\geq D$. Then $C_1+iA \equiv C_2+iA$.
\end{lem}

\noindent{\sl Proof:} By Proposition
\ref{prop:A+iB=C+iD}
we have  
\be\label{eq:C-stratum-cond}
\Hc^{C_1} \cap\Hc^D \,=\,\Hc^{C_2}\cap \Hc^D. 
\ee
If  $A\geq D$, then $\Hc^A\subset\Hc^D$, so intersecting \eqref{eq:C-stratum-cond}
with $\Hc^A$, we get $\Hc^{C_1}\cap\Hc^A=\Hc^{C_2}\cap \Hc^A$, i.e., that
$C_1+iA \equiv C_2+iA$. 
This proves Lemma \ref{lem:C-stratum-enlarge} and the condition (Q). 

\vskip .2cm

We now prove that $\Ec^\bullet(E)$ is $\tau\Sc^{(1)}$-constructible. For this it suffices to prove that
the shifted Verdier dual $\Ec^\bullet(E)^*$ is $\tau\Sc^{(1)}$-constructible. But writing $\Ec^\bullet(E)$
as the total object of the double complex as in Corollary \ref{cor:E(F)-bicell}, and applying the duality
term by term, we find that $\Ec^\bullet(E)^*$ is the total object of a double complex as in 
\eqref{double-resol-dual} which is the same as the complex of sheaves $\wc\Ec^\bullet(E^*)$
corresponding to the dual matrix diagram $E^*$ and the system of tube cells $B+i\RR^n$ obtained from the
 system of the $\RR^n+iA$ by applying $\tau$.  So 
 \be\label{eq:EE*=tEE*}
 \Ec^\bullet(E)^* \,\simeq \, \wc\Ec^\bullet(E^*) \,=\, \tau^{-1} \Ec^\bullet(E^*)
 \ee
  is   $\tau \Sc^{(1)}$-constructible because $\Ec^\bullet(E^*)$ is $\Sc^{(1)}$-constructible by the above. 
 This finishes the proof of Proposition \ref{prop:E(E)S0constr}. 
 
 \begin{prop}
 The complex $\Ec^\bullet(E)$ is perverse. We have therefore a functor
 \[
 \GG: \Mc_\Hc\lra\Perv(\CC^n,\Hc), \quad E\mapsto \Ec^\bullet(E). 
 \]
 \end{prop}
 
 \noindent{\sl Proof:} We first prove the condition $(P^-)$ of perversity: that $\ul H^p(\Ec^\bullet(E))$
 is supported on a complex submanifold of complex codimension $\geq p$. By construction, the $p$th term
 $\Ec^p(E)=\bigoplus_{\cod(A)=p} \Ec_A(E)$ is supported on the union of the $\RR^n+ \ol A$ 
 where $A$ runs over faces of $\Hc$ of real codimension $p$. So $\Ec^p(E)$ and therefore $\ul H^p(\Ec^\bullet(E))$
 is supported on the union of $\RR^n+iL$ where $L$ runs over flats of $\Hc$ of real codimension $p$.
 But since, by Proposition \ref{prop:E(E)S0constr}, $\Ec^\bullet(E)$ is $\Sc^{(0)}$-constructible,
 $\Supp \, \ul H^p(\Ec^\bullet(A))$ is a complex manifold, in fact, a finite-union of $\CC$-linear subspaces
 $M\subset\CC^n$. But if such an  $M$ lies in $\RR^n+iL$, it must lie in $L+iL=L_\CC$ which has complex codimension $p$. 
 This proves  $(P^-)$ for $\Ec^\bullet(E)$.
 
 Now,  $(P^+)$ is equivalent to  $(P^-)$  for $\Ec^\bullet(E)^*$. By \eqref {eq:EE*=tEE*}, we have
 $ \Ec^\bullet(E)^*= \tau^{-1}\Ec(E^*)$, and $(P^-)$   for $\Ec(E^*)$ has just been proved.
 \qed
 
 \begin{prop}
 The functors $\EE$ and $\GG$ are quasi-inverse to each other and so are equvalences of categories.
 \end{prop}
 
 \noindent{\sl Proof:} 
 That $\GG\circ \EE\simeq\Id$ is clear: the Cousin complex of $\Fc$ is a representatve of $\Fc$.
 Conversely, suppose we start from a matrix diagram $E=(E_{A,B}, \del', \del'')$
 and form the   complex $\Ec^\bullet=\Ec^\bullet(E)$ whose data, as a complex of cellular sheaves, is completely
 equivalent to the bicomplex as in Corollary \ref{cor:E(F)-bicell}, i.e., yo $E$. We need to
 prove that the ``intrinsic Cousin complex" associated to $\Ec^\bullet$, is $\Ec^\bullet$ itself. This argument is
 elementary and similar to \cite{KS}, \S 6. 
 
 More precisely, for a face $D$ let $\lambda_D: \RR^n+iD\hra\CC^n$ be, as before,  the embedding.
 It is enough to prove that for any $k$-vector space $V$ (we will need $V=E_{AB}$)
 \[
 \lambda_{D*}\, \lambda_D^! \bigl( \eps^{BA}_{!*} \, \ul{V}_{B+iA}\bigr) \,=\,
 \begin{cases}
 0, &\text{ if } D\neq A,
 \\
 \eps^{BD}_{!*}\,  \ul{V}_{B+iD}, & \text {if } D=A. 
 \end{cases}
 \]
 This reduces to the case $V=\k$ which is a Cartesian product situation. So we reduce to a statement about
 the second factor only, that is, denoting by $j_C: C\to\RR^n$ the embedding of a face $C\in\Cc$, that
 \[
 j_D^!\,  j_{A*} \, \ul\k_A \,=\,0, \quad \text{if } D\neq A. 
 \]
(of course, the LHS is equal to $\ul \k_D$, of $D=A$). But this claim is clear: by Verdier duality, it is equivalent to
\[
j_D^*\,  j_{A!} \, \ul\k_A =0, \quad \text{if } D\neq A,
\] 
which is completely obvious, as $j_{A!}\, \ul\k_A$ is just the extension of the constant sheaf by $0$ from $A$ to
$\ol A$ and then to $\RR^n$. \qed

\vskip .2cm

This finishes the proof of Theorem \ref{thm:main}.


 \section{Examples and complements}\label{sec:ex-comp}
 
 \paragraph{The $1$-dimensional case.} 
 Let $n=1$ and let $\Hc$ consist of the ``hyperplane'' $0\in\RR$. 
 The corresponding category $\Perv(\CC,0)$ consists of perverse sheaves on $\CC$
 with the  possible singularity at $0$. 
 
 The poset $\Cc$ of faces has $3$ elements:  $\RR_-, \{0\}$ and $\RR_+$, so a matrix diagram
 has the form
 \eqref {eq:perv-C-0} or, with the notations for the arrows spelled out,
 \be\label{eq:perv-C-0-2}
\xymatrix{
E_{-,+} \ar[r]^\simeq_{\del''}& E_{0,+}& \ar[l]_{\simeq}^{\del''} E_{+,+}
\\
E_{-,0}\ar[r]^{\del''}
\ar[u]^{\simeq}_{\del'}  \ar[d]_{\simeq}^{\del'}
& E_{0,0}
\ar[u]^{\del'} \ar[d]_{\del'}
 & \ar[l]_{\del''} E_{+,0}
 \ar[u]_\simeq^{\del'}  \ar[d]^\simeq_{\del'}
\\
E_{-,-}\ar[r]_\simeq^{\del''}  &E_{0,-} & \ar[l]^{\simeq}_{\del''}  E_{+,-}
}
\ee
  Theorem \ref{thm:main} says therefore that
 $\Perv(\CC,0)$ is equivalent to the category of diagrams  \eqref {eq:perv-C-0-2}. 
 Let us compare this with other known descriptioms. The most classical description
 \cite{beil-gluing, galligo-GM} is  in terms of diagrams of vector spaces
   \be\label{eq:phi-psi}
  \xymatrix{
   \Phi \ar@<.7ex>[r]^a& \Psi \ar@<.7ex>[l]^b
  }, \quad \text {  $\on{Id}_\Psi-ab$ is invertible.}
    \ee
  The approach of  \cite{KS}  gives rise to
   another description, in terms of diagrams of vector spaces 
\be\label{eq:dirac}
\begin{gathered}
 \xymatrix{
E_-\ar@<-.7ex>[r]_{\delta_-}& E_0
\ar@<-.7ex>[l]_{\gamma_-}
\ar@<.7ex>[r]^{\gamma_+}& E_+
\ar@<.7ex>[l]^{\delta_+}
}
 \\
\gamma_-\delta_- = \on{Id}_{E_-}, \,\,\,\gamma_+
\delta_+ = \on{Id}_{E_+},\\
\gamma_-\delta_+: E_+\to E_-, \,\,\, \gamma_+\delta_-: E_-\to E_+ \quad \text {are invertible}. 
\end{gathered}
\ee
 See \cite{KS} \S 9A for a direct constriction of an equivalence between the categories of  diagrams
 \eqref{eq:phi-psi} and \eqref{eq:dirac}.  Let us explain an equivalence between the
 categories of diagrams \eqref {eq:perv-C-0-2} and  \eqref{eq:dirac}. 
 Given a  diagram as in  \eqref {eq:perv-C-0-2}, we consider first its
 middle horizontal part (the $0$th row) which gives the straight arrows  $\delta_+,\delta_-$ in
 \be\label{eq:dirac-curved}
 \xymatrix{
 E_{-,0} \ar[r]_{\delta_- = \del''} & E_{0,0}
 \ar@/_2pc/[l]_{\gamma_-}
 \ar@/_2pc/[r]_{\gamma_+}
 & E_{+,0} \ar[l] _{\delta_+=\del''}
 }
 \ee
 The curved  arrows $\gamma_\pm$  are defined  as the compositions along the
 corresponding squares in  \eqref {eq:perv-C-0-2}. That is, $\gamma_-$
 is the composition 
 \[
 E_{0,0}\buildrel\del'\over \lra E_{0,+} \buildrel(\del'')^{-1}\over  \lra E_{-,+}
 \buildrel(\del')^{-1}\over \lra E_{-,0},
 \]  
 while $\gamma_+$ is the composition
 \[
 E_{0,0}\buildrel \del'\over\lra E_{0,-}
 \buildrel (\del'')^{-1}\over\lra E_{+,-}\buildrel (\del')^{-1}\over\lra E_{+,0}. 
 \]
 The commutativity of  \eqref {eq:perv-C-0-2} and invertibility of
 the arrows at its outer rim implies easily  that the diagram \eqref{eq:dirac-curved}
 is of the type \eqref{eq:dirac}. Further,   this procedure
 gives an equvialence between the categories of diagrams  \eqref {eq:perv-C-0-2}
 and  \eqref{eq:dirac}. We leave the verifications to the reader. 
 
 \paragraph{Comparison with \cite{KS}.} 
 For a general arrangement $\Hc$ of linear hyperlplanes in $\RR^n$ we gave, in
 \cite{KS}, a description of $\Perv(\CC^n,\Hc)$ in terms of ``single-indexed'' diagrams
 \[
Q\,=\,  \bigl( (E_A)_{A\in\Cc}, (\gamma_{A_1,A_2}, \delta_{A_2, A_1})_{A_1\leq A_2}\bigr),
 \]
 where:
 \begin{itemize}
 \item $E_A$ are finite-dimensional $\k$-vector spaces given for any $A\in\Cc$.
 
 \item $\gamma_{A_1, A_2}: E_{A_1}\to E_{A_2}$, resp.
 $\delta_{A_2, A_1}: E_{A_2}\to E_{A_1}$ are linear maps forming a representation,
 resp. an anti-representation of $\Cc$ on  $(E_A)$ and satisfying
 the axioms of monotonicity ($\gamma_{A_1, A_2}\delta_{A_2, A_1}=\Id$), 
transitivity and invertibility, see \cite{KS}.  
 \end{itemize}
Let us compare this with the ``matrix'' description given by Theorem \ref{thm:main}. 
By the construction of \cite{KS}, the partial data $(E_A, \gamma_{A_1, A_2})$
are just the linear algebra data describing the cellular sheaf
$i_\RR^!\Fc[n]$ on $\RR^n$, where $i_\RR:\RR^n\hra\CC^n$ is the embedding. Therefore
\[
E_A = E_{0,A}, \,\,\,\gamma_{A_1, A_2}=\del''
\]
is just one column of the matrix diagram $(E_{A,B},\del', \del'')$ corresponding to $\Fc$. 

Further, any $E_{A,B}$ is, by our construction, the stalk at $B+iA$ of the sheaf
$j_{A*} j_A^!\otimes\orr(A)[\cod(A)]$, where $j_A: \RR^n+iA\hra \CC^n$ is the embedding.
There stalks were  described in \cite{KS}, Cor. 3.22, and we get
\be\label{eq:EAB=EBcA}
E_{A,B} \,=\, E_{B\circ A},
\ee
where $B\circ A$ is ``the first cell in the direction $A$ visible from $B$. 
More precisely, (see \cite{KS} Prop. 2.3) $B\circ A$ is the (uniquely defined) cell containing the points
\be\label{eq:circ}
(1-\eps)b  +  \eps a, \quad b\in B,\, a\in A,\,\, 0<\eps \ll 1. 
\ee

In particular, $E_{A,0}$ is also identified with $E_A$, and the maps $\del'$ connecting
different $E_{A,0}$, are precisely the $\delta_{A_2, A_1}$, as both appear from the
differential in the Cousin complex. So the partial data $(E_A,\delta_{A_2, A_1})$
is just the $0$th row of  $(E_{A,B})$.

\vskip .2cm

The operation $\circ$ was introduced by Tits \cite{tits} in 1974 in the context of buildings
 (which includes arrangements of root hyperplanes) and later, independently, by 
Bj\"orner, Las Vergnas, Sturmfels,
White and Ziegler \cite{BLSWZ} for oriented matroids (which includes all hyperplane arrangements).
For simplicity, we will refer to $\circ$ as the {\em Tits product}. Let us list some of its properties.

\begin{prop}\label{prop:circ}
\begin{itemize}
\item[(a)] The Tits product $\circ$ is associative (but not commutative).

\item[(b)] Further,  $\circ$ is monotone in the second argument (but not in the first one). That is, if
$A_1\leq A_2$, then $B\circ A_1\leq B\circ A_2$ for any $B$. 

\item[(c)] In the notation of \eqref{eq:H^C} we have
\[
\Hc^{B\circ A} \,=\, \Hc^B \cap \Hc^A. 
\]
In particular (by Proposition \ref{prop:A+iB=C+iD}), $B\circ A$ and $A\circ B$ lie in the
same complex stratum.  

\item[(d)] Let $B, A_1, A_2\in\Cc$ be such that $A_1\leq A_2$. Then the following are equivalent:
\begin{itemize}
\item[(di)]  $B+iA_1$ and $B+iA_2$ lie in the same complex stratum.
\item[(dii)] $B\circ A_1 = B\circ A_2$. 
\end{itemize}

\end{itemize}
\end{prop}

 \noindent{\sl Proof:} Parts (a) and (b) are well known, see, e.g., \cite{KS} Props. 2.3(a) and
 2.7(a).  Part (c) follows at once from \eqref{eq:circ}.  Let us prove (d). Suppose
 (di) holds. Then, by Proposition \ref{prop:A+iB=C+iD},
 \be\label{eq:BA1=BA2}
 \Hc^B\cap \Hc^{A_1} \,=\, \Hc^B\cap\Hc^{A_2}. 
 \ee
Since $A_1\leq A_2$, we have $B\circ A_1\leq B\circ A_2$ by (b). But \eqref{eq:BA1=BA2}
means, in virtue of (c), that $\Hc^{B\circ A_1} = \Hc^{B\circ A_2}$, and this implies
that $\dim (B\circ A_1)=\dim(B\circ A_2)$. Therefore we must have $B\circ A_1 = B\circ A_2$,
that is, (dii) holds. 

Conversely, suppose (dii) holds. Then, by (c), we have \eqref{eq:BA1=BA2} so
$B+iA_1$ and $B+iA_2$ lie in the same complex stratum by Proposition 
 \ref{prop:A+iB=C+iD}, that is, (di) holds.\qed

 \begin{rems}\label{rems:circle-cond}
 
 (a) Part (c) of  Proposition  \ref{prop:circ} can be compared, via \eqref{eq:EAB=EBcA}, with Remark \ref{rem:A+ib=B+iA}. That is, even though $B\circ A\neq A\circ B$ in general,
 $E_{B\circ A}$ is isomorphic to $E_{A\circ B}$. 
 
 \vskip .2cm
(b)  Part (d) of Proposition \ref{prop:circ} means that the conditions in the 
axioms ($M3'$) and ($M3''$)
 of a matrix diagram can be formulated in terms of the Tits product. More precisely, the conditions that
 $B+iA_1 \equiv B+iA_2$ (i.e., $B+iA_1$ and $B+iA_2$ lie in the same complex stratum) in
  in ($M3''$) directly means that $B\circ A_1=B\circ A_2$. The condition that $B_1+iA \equiv B_2+iA$ 
  in ($M3'$) is equivalent to $A+iB_1 \equiv A+iB_2$, i.e., to $A\circ B_1 = A\circ B_2$. 
 
 \end{rems}

To summarize,  the matrix diagram $(E_{A,B})$ contains the same vector spaces as the
single-indexed one $(E_A)$ but with repetitions, being a kind of ``Hankel matrix''
with respect to the  Tits product $\circ$. It is these repetitions that allow us to
write the relations among the arrows $\del', \del''$ of a matrix diagram in such a
simple, local form: as commutativity of elementary squares. 

\paragraph{Affine arrangements.}
Let now $\Hc$ be a, possibly infinite,  arrangement of {\em affine hyperplanes} in $\RR^n$. 
For any affine hyperplane $H\in\Hc$ with real affine equation $f_H(x)=a$, where
$f_H: \RR^n\to\RR$ is $\RR$-linear nd $a\in\RR$, let $\ol H\subset\RR^n$ be the linear hyperplane
with the equation $f_H(x)=0$. We denote $\ol\Hc$ the linear arrangement of the
hyperplanes $\ol H$, $H\in\Hc$ (ignoring possible repetitions) and assume that:
\begin{itemize}
\item $\Hc$ is  closed (as a subset in $\RR^n$) and  locally finite, i.e., any $x\in\RR^n$
has a neighborhood  meeting only finitely many affine hyperplanes from $\Hc$.

\item $\ol\Hc$ is finite. 
\end{itemize}

 The concepts of flats of $\Hc$, their complexifications and the stratification $\Sc^{(0)}$
 of $\CC^n$ into generic parts of complexified flats are defined analogously to the
 case of linear arrangements. We then have the category $\Perv(\CC^n,\Hc)$ of
 perverse sheaves on $\CC^n$ smooth with respect to $\Sc^{(0)}$. Let us
 give a modification of Theorem \ref{thm:main} to the case of affine arrangements as above.
 
 \vskip .2cm
 
 We denote $\Sc^{(2)}$ the quasi-regular cell decompostion of $\CC^n$ into product cells of the form
 $iA+B$ with $A\in\ol\Cc$ and $B\in\Cc$. 
 
 \begin{prop}
 The decomposition $\Sc^{(2)}$ refines $\Sc^{(0)}$.
 \end{prop}
 
\noindent{\sl Proof:} This is a consequence of the following obvious remark.
Let $f: \RR^n\to\RR$ be an $\RR$-linear function and $a\in\RR$. Denote
by $f^\CC:\CC^n\to\CC$ the complexification of $f$. Then, for $x,y\in\RR^n$ the
condition $f^\CC(x+iy)=a$ is equivalent to $f(x)=a$ amd $f(y)=0$. \qed

\begin{thm}\label{thm:affine}
The category $\Perv(\CC^n,\Hc)$ is equivalent to the category of
 diagrams of finite-dimensional $\k$-vector spaces of the form
 \[
 \begin{gathered}
 E_{A,B}\in\Vect_\k, \,\, A\in \ol\Cc, \, B\in \Cc,
 \\
  \del':_{(A|B_1, B_2)}:  E_{A, B_1}\lra E_{A, B_2},\,\, B_1\leq B_2,
 \\
  \del''_{(A_1, A_2|B)} : E_{A_1,B}\lra E_{A_2,B},\,\, A_1\leq A_2,
  \end{gathered}
 \]
such that:
\begin{itemize}
\item[(A2)] The maps $\del',\del''$ define a representation of $\ol\Cc^\op\times\ol\Cc$ in $\Vect_\k$. 

\item[(A3)] If $iA+B_1$ and $iA+B_2$ lie in the same stratum of $\Sc^{(0)}$, then
$\del'_{(A|B_1, B_2)}$ is an isomorphism. If $iA_1+B$ and $iA_2+B$ lie in the same stratum of
$\Sc^{(0)}$, then $\del''_{(A_1, A_2|B)}$ is an isomorphism. 
\end{itemize}
\end{thm}

\noindent {\sl Proof:} It can be obtained, as in \cite{KS} \S 9B, by an amalgamation argument from the linear case, using
the fact that perverse sheaves form a stack of categories. Alternatively, one can perform
a direct analysis of the
Cousin complex associated to $\Fc\in\Perv(\CC^n,\Hc)$ and formed by the sheaves
\[
j_{A*} j_A^! \Fc[\cod(A)], \,\, A\in\ol\Cc, \quad j_A: \RR^n+iA \hra\CC^n. 
\]
\qed

\begin{ex}
Consider the arrangement of two points $0,1$ in $\RR$. Then
\[
\begin{gathered}
\Cc \,=\, \{ \RR_{<0}, 0, (0,1), 1, \RR_{>1}\}, \quad \ol\Cc \,=\, \{\RR_-, 0, \RR_+\}, 
\end{gathered}
\]
The decomposition $\Sc^{(2)}$  of $\CC$ into the product cells is depicted in Fig. 
\ref{fig.affine}.

 \begin{figure}[H]
 \centering
 \begin{tikzpicture}[scale=.4, baseline=(current bounding box.center)]
 
 \node at (0,0){$\bullet$};
 \node at (4,0){$\bullet$}; 
 \draw[line width=0.9] (-4,0) -- (10,0); 
 \draw (0,-4) -- (0,4); 
  \draw (4,-4) -- (4,4); 
  \node at (-0.5,-0.5){$0$}; 
  \node at (11,0){$\RR$}; 
  \node at (2, 0.5){$(0,1)$}; 
  \node at (7,0.7){$\RR_{>1}$}; 
  \node at (-3,0.7){$\RR_{<0}$}; 
 
   \node at (3.5,-0.5){$1$}; 
  \end{tikzpicture}
 \caption{ Product cells for the arrangement $\Hc=\{0,1\}\subset\RR$.}\label{fig.affine}
 \end{figure}
 
 Theorem \ref{thm:affine} identifies $\Perv(\CC,\Hc)$ with the category of commutative diagrams 
 of the form
 \[
 \xymatrix{
 E_{<0,+} \ar[r]^\simeq & E_{0,+} & \ar[l]_\simeq E_{(0,1), +} \ar[r]^\simeq & E_{1,+} &
 \ar[l]_\simeq E_{>1, +}
 \\
  E_{<0,0}\ar[u]^\simeq  \ar[d]_\simeq  \ar[r] & E_{0,0}\ar[u] \ar[d]  & \ar[l]  E_{(0,1), 0} \ar[u]_\simeq 
  \ar[d]^\simeq
  \ar[r]  & 
  E_{1,0}   \ar[u] \ar[d] &
 \ar[l] E_{>1, +}\ar[u]_\simeq \ar[d]^\simeq
 \\
  E_{<0,-} \ar[r]_\simeq & E_{0,-} & \ar[l]^\simeq E_{(0,1), -} \ar[r]_\simeq & E_{1,-} &
 \ar[l]^\simeq E_{>1, -}
}
 \]
 Such a diagram can be seen as an amalgamation of two diagrams of the form 
 \eqref{eq:perv-C-0-2}. 
 \end{ex}


\vskip 0.1cm

M.K.: Kavli IPMU, 5-1-5 Kashiwanoha, Kashiwa, Chiba, 277-8583 Japan,
{\tt mikhail.kapranov@ipmu.jp}

V.S.: Institut de Math\'ematiques de Toulouse, Universit\'e Paul Sabatier, 118 route de Narbonne, 
31062 Toulouse, France, 
 {\tt schechtman@math.ups-tlse.fr }

\end{document}